\documentclass{article}
\usepackage[utf8]{inputenc}

\usepackage[english]{babel}

\usepackage[utf8]{inputenc}
\usepackage[T1]{fontenc}

\usepackage{amssymb}
\usepackage{stmaryrd}
\usepackage{amsmath,amsfonts}
\usepackage[tikz]{bclogo}
\usepackage[top=4.34cm]{geometry}
\usepackage{subcaption}

\usepackage{lmodern}
\usepackage{textcomp}
\usepackage{mathtools, bm}
\usepackage{amssymb, bm}
\usepackage[unq]{unique}
\usepackage{enumerate}
\usepackage{comment}
\usepackage[breaklinks]{hyperref}

\usepackage[numbers]{natbib}  

\usepackage[breaklinks]{hyperref}
\usepackage[numbers]{natbib}

\usepackage{amsthm}

\usepackage{cleveref}

\newtheorem{theorem}{Theorem}[section]
\newtheorem{lemma}[theorem]{Lemma}

\newtheorem{conjecture}[theorem]{Conjecture}

\newtheorem{Question}[theorem]{Question}
\newtheorem{Problem}[theorem]{Problem}
\newtheorem{Claim}[theorem]{Claim}

\theoremstyle{definition}

\newcommand{\cC}{\mathcal{C}}

\title{A proof of the 3/4 conjecture for the total domination game}

\author{Julien Portier\footnote{\href{mailto:jp899@cam.ac.uk}{jp899@cam.ac.uk}, Department of Pure Mathematics and Mathematical Statistics (DPMMS), University of Cambridge, Wilberforce Road, Cambridge, CB3 0WA, United Kingdom} \and Leo Versteegen \footnote{\href{mailto:lvv23@cam.ac.uk}{lvv23@cam.ac.uk}, Department of Pure Mathematics and Mathematical Statistics (DPMMS), University of Cambridge, Wilberforce Road, Cambridge, CB3 0WA, United Kingdom}}

\date{}

\begin{document}

\maketitle

\begin{abstract}
    The total domination game is a game played on a graph $G$ by 2 players, Dominator and Staller, who alternate in selecting vertices until each vertex in the graph $G$ has a neighbor in the set of selected vertices. Dominator's aim is to arrive at this state in as few moves as possible, while Staller wants to achieve the opposite. In this paper we prove the 3/4-conjecture by Henning, Klav{\v{z}}ar, and Rall, stating that Dominator has a strategy to finish the total domination game on $G$ within $3n/4$ moves, where $n$ is the number of vertices of $G$. We also prove that if the minimum degree of $G$ is at least 2, then Dominator can finish the game within $5(n+1)/7$ moves.
\end{abstract}

\section{Introduction}

Let $G$ be a graph without isolated vertices. A vertex is \emph{totally dominated} by a set $A\subset V(G)$, if it has a neighbor in $A$. A totally dominating set of $G$ is a set that totally dominates every vertex of $G$. The total domination game was introduced as a natural variant of the domination game by Henning, Klav{\v{z}}ar and Rall in \cite{henning2015total}. The two players, Dominator (D) and Staller (S), alternately select vertices of a graph $G$ without isolated vertices. Each selected vertex must strictly increase the number of vertices that are totally dominated by the set of selected vertices. The game stops when the set of selected vertices is a total dominating set of $G$. Dominator's aim is to minimize the number of vertices in this set at the end of the game, while Staller's aim is to maximize it. The \emph{game total domination number} $\gamma_{tg}(G)$ of $G$ is defined as the number of vertices in the resulting set when Dominator starts the game and both players play optimally. In \cite{henning2017_0.8}, Henning, Klav{\v{z}}ar, and Rall proved that $\gamma_{tg}(G) \leq \frac{4}{5}n$ for every graph $G$ without isolated vertices or edges, and they conjectured that the following holds.

\begin{conjecture}\label{3-4_conj}
For every graph $G$ without isolated vertices or edges, $\gamma_{tg}(G) \leq \frac{3}{4}n$.
\end{conjecture}

\Cref{3-4_conj} is plainly called the \emph{3/4-conjecture}, and it has been verified for graphs of minimum degree at least 2 in \cite{bujtas2016transversal} and for graphs with the properties that $\deg(u) + \deg(v) \geq 4$ for every edge $uv \in E(G)$ and that no two vertices of degree one are at distance four in \cite{henning2016progress}. The best bound for arbitrary graphs without isolated vertices or edges so far was due to Bujt{\'a}s \cite{bujtas2018game}, who proved for every such graph the bound $\gamma_{tg}(G) \leq \frac{11}{14}n$. In this paper we prove \Cref{3-4_conj}.

\begin{theorem}\label{3-4_thm}
For every graph $G$ without isolated vertices or edges, $\gamma_{tg}(G) \leq \frac{3}{4}n$.
\end{theorem}

Note that \Cref{3-4_thm} is best possible, as the example of the union of disjoint copies of $P_4$ shows. The admission of isolated edges also adds nothing of interest since both vertices of such edges will ultimately have to be played. \\

As mentioned before, in \cite{bujtas2016transversal}, Bujt{\'a}s, Henning and Tuza verified \Cref{3-4_conj} for graphs of minimum degree at least 2, and in fact, they proved that for every such graph $G$, $\gamma_{tg}(G) < \frac{8}{11}n$. In this paper, we improve on their result as follows.

\begin{theorem}\label{Bound_Min_degree_2}
For every graph $G$ with minimum degree at least 2, $\gamma_{tg}(G) \leq \frac{5}{7}n+\frac{5}{7}$.
\end{theorem}

\section{Preliminaries}

\subsection{Main idea of the proofs}

Before we go into the details of the proofs of \Cref{Bound_Min_degree_2} in Section 3 and of \Cref{3-4_thm} in Section 4, we want to convey their central idea. For a graph $G$ on $n$ vertices, we can make the following two observations.
\begin{enumerate}
    \item At any point in the game, the number of totally dominated vertices is at most $n$.
    \item At any point in the game, the sum of the number of played vertices and the number of vertices that have not been played is $n$.
\end{enumerate}

Moreover, if at some point in the game, a vertex $v$ has not been played and the last neighbor of $v$ becomes totally dominated, then $v$ is \emph{unplayable} from there on out. To exploit the two observations above, we would like Dominator to play in a way such that the sum of the numbers of totally dominated and unplayable vertices grows as quickly as possible. Whether it is easier to totally dominate many new vertices or to make many of them unplayable, changes over the course of the game. For this reason, the strategies of D consist of a number of different \emph{phases} that follow each other linearly and come with separate sets of instructions for D. In each phase, we will claim and prove bounds on the number of new totally dominated and unplayable vertices depending on the number of moves that were played during the phase that D can achieve by playing according to our instructions. At the end of the proofs of \Cref{3-4_thm} and \Cref{Bound_Min_degree_2}, we will then combine these bounds to show that $\gamma_{tg}(G)\leq 3n/4$ in the general case and $\gamma_{tg}(G)\leq 5(n+1)/7$ in the minimum degree at least $2$ case, respectively.

\subsection{Describing the state of the game}

Intuitively, S should always wish that as few vertices as possible get totally dominated, but suboptimal moves played by S or unexpectedly good opportunities for D can often lead to more complicated case distinctions later on. For this reason, we will sometimes want to pretend that certain vertices are not yet dominated when they in fact are. To leave no doubt about the validity of such arguments, we substitute some natural notions of the game, such as a vertex being totally dominated by artificial counterparts that track how often we explicitly registered these things during the game. These artificial concepts are not part of the game, but rather auxiliary pieces of data kept by us as observers of the game.

Specifically, we will think of all vertices as \emph{painted} in one of the colors \emph{white} or \emph{black}, and in addition, we may mark them as \emph{depleted} or \emph{dependent}. We impose on ourselves the following restrictions with respect to how vertices may be painted or marked.

\begin{itemize}
    \item If a vertex is not totally dominated, it has to be white.
    \item Once painted black, a vertex is not allowed to be re-painted in white.
    \item A vertex is only allowed to be depleted if it has not been played and it has no white neighbors.
    \item If a vertex has been marked as depleted or dependent, we are not allowed to remove this mark.
    \item A vertex is only allowed to be dependent if it has at most one white neighbor.
    \item A vertex is not allowed to be marked as both depleted and dependent.
\end{itemize}

In particular, every vertex is painted white at the start of the game, and at any stage of the game, the number of black vertices is always bounded from above by the number of totally dominated vertices. Note further that once a vertex has been marked as depleted, it cannot be played after that.

As usual, a \emph{leaf} is a vertex of degree exactly 1. The \emph{parent} of a leaf is its unique neighbor. A point $t$ of the game refers to a situation in our analysis after $t$ moves have been played, always assuming that D has played according to our instructions thus far. In the proofs of both \Cref{3-4_thm} and \Cref{Bound_Min_degree_2}, we wish to keep track of the number of black vertices $\beta(t)$ and the number of depleted vertices $\delta(t)$ at any given point in time $t$. For the general case of graphs with leaves, we will also keep track of the number of leaves that have been played or marked as depleted, which we denote by $\lambda(t)$. By the rules we have set about black, depleted and dependent vertices, it is easy to see that $\beta,\lambda$ and $\delta$ are increasing in $t$.

As for dependent vertices, we denote their number by $\sigma$ and the number of undominated neighbors of dependent vertices by $\nu$. For the proof of \Cref{3-4_thm}, we are only interested in their difference $\chi=\sigma-\nu$. To prove the stronger bound in \Cref{Bound_Min_degree_2} however, we will also involve the quantity $\nu$ directly. We point out that $\chi$ too is monotonous.

\section{Proof of \Cref{Bound_Min_degree_2}}

The strategy is divided into consecutive \emph{phases}, each of them starting after the last one ends by meeting a predefined end condition. We track the number of moves played so far by a variable $t$. There are going to be four phases, and for each $i\in [4]$, $t_i$ is the last move before the end of Phase $i$ and let $t_0=0$. It can happen that phases are skipped because their end conditions are met from the start of the phase and in this case we have $t_i=t_{i-1}$. For each $i\in [4]$, let $T_i=t_i-t_{i-1}$.

\subsection{Phase 1}\label{SectionPhase1}
This phase starts at the start of the game and ends once it is Dominator's turn and no unplayed vertex has at least three white neighbors. For a set $A\subset V$, a \emph{triple separation} is a family $(B_v)_{v\in A}$ such that for all $v\in A$, $B_v\in N(v)^{(3)}$ and for all distinct $v,w\in A$, $B_v\cap B_w=\emptyset$. The triple separation is \emph{white}, if for all $v\in A$ all elements of $B_v$ are white, and a set $A$ is called \emph{white-separable} if no vertex in $A$ has been played and $A$ has a white triple separation. Let $M$ denote the size of the largest set $A$ for which there exists a triple separation.

\begin{Claim}\label{phase:3+}
D can play in a way such that the following hold.
\begin{itemize}
    \item No vertex is marked as dependent.
    \item $\beta(t_1) \geq 2T_1$.
    \item $T_1\geq M$.
\end{itemize}
\end{Claim}
\begin{proof}
Let $M_t$ be the size of the largest white-separable set after move $t$. By definition of $M$ and due to the fact that all vertices are white at the start of the game, $M_0=M$. The phase is over exactly when $M_t=0$ and we will ensure that the phase lasts for at least $M$ moves by maintaing that $M_t\geq M-t$. Whenever it is Dominator's turn, D plays a vertex $v$ from a largest white-separable set $A$, and fixing some white triple-separation $(B_v)_{v\in A}$, we paint the vertices in $B_v$ black. Whenever it is Staller's turn and S plays a vertex $v$, we paint one arbitrarily chosen neighbor of $v$ black. In both cases, we have $M_t\geq M_{t-1}-1$.

Note that after every two moves four new vertices are painted black. Therefore, we have $\beta(t_1)\geq 2T_1$.
\end{proof}

\subsection{Phase 2}\label{sec:phase2-mindeg2}
For $k\in \mathbb{N}$, a \emph{$k$-walk} is a sequence of vertices $w_1v_1\ldots w_{k}v_{k}$ such that for all $1\leq i<j\leq k$, $v_i\neq v_j$, $w_i\neq w_j$, $w_i$ is white, $v_i$ is adjacent to $w_{i+1}$ and $w_j$ is adjacent to $v_j$. A $k$-walk is a \emph{$k$-circuit} if $v_k$ is adjacent to $w_1$ and a walk is \emph{terminal} if it is a $2$-walk and $v_2$ has no other white neighbor than $w_2$.

The current phase ends once it is Dominator's turn, and no 6-walk, terminal walk or 2-circuit exists. Note that after this phase, no 6-walk or 2-circuit can be reintroduced.

\begin{Claim}\label{phase:long-walks}
D can play in a way such that the following hold.

\begin{itemize}
    \item $\beta(t_2)-\beta(t_1)\geq 3T_2/2$.
    \item $4\delta(t_2)\geq T_2$.
    \item No vertex is marked as dependent.
\end{itemize}
\end{Claim}

\begin{proof}
We show that if at some point $t$ it is Dominator's turn and the phase is not over, at least one of the following hold.
\begin{itemize}
    \item D can play the next move so that
    \begin{align*}
        \beta(t+2)-\beta(t)\geq 3 \qquad \text{and} \qquad \delta(t+2)-\delta(t)\geq 1,
    \end{align*}
    \item D can play the next two moves so that 
    \begin{align*}
    \beta(t+4)-\beta(t)\geq 6 \qquad \text{and}\qquad
    \delta(t+4)-\delta(t)\geq 1
    \end{align*}
\end{itemize}
Once this dichotomy is shown, a simple inductive argument shows the Claim. Before the first move of its phase, and after each of these sequences, we paint all vertices that are white and totally dominated black.

Suppose first there exists a 6-walk $w_1\ldots v_6$. Then D plays $v_1$. If in the subsequent move, S plays some vertex $x$ that is not adjacent to $w_3$ or $w_4$, D responds by playing $v_3$ after which S plays some second move $y$. We then paint $w_1,w_2,w_3,w_4$ and all white neighbors of $x$ and $y$ black and mark $v_2$ as depleted. Similarly, if S responds by playing a vertex that is adjacent to $w_3$, D plays $v_4$ next and we can mark $v_3$ as depleted. Finally, if S plays a vertex adjacent to $w_4$, D can play $v_5$ next. 

Suppose next there exists a 2-circuit $w_1v_1w_2v_2w_1$. D simply plays $v_1$, we paint $w_1$ and $w_2$ black and we mark $v_2$ as depleted. Whichever vertex Staller plays, it is now sufficient to paint all its neighbors black to arrive at the first outcome.

Suppose finally that there exists a terminal walk $w_1v_1w_2v_2$. Here, D plays $v_1$, we color $w_1$ and $w_2$ black, and we mark $v_2$ as depleted. Once again, whichever vertex Staller plays, we arrive at the first outcome.
\end{proof}

\subsection{Phase 3}\label{sec:phase3-mindeg2}
This phase ends once it is Dominator's turn, no terminal walk exists, no white vertex with degree higher than two and part of a 2-walk exists and no 5-walk exists. Recall that $M$ is the largest set with a triple separation.

\begin{Claim}\label{phase:intersecting-circuits}
D can play in a way such that the following hold.

\begin{itemize}
    \item $\beta(t_3)-\beta(t_2)\geq 3T_3/2$.
    \item $4(\delta(t_3)-\delta(t_2))+\chi(t_3)\geq T_3$
    \item $\nu(t_3)\leq (\chi(t_3)+M)/2$.
    \item All dominated vertices are black.
\end{itemize}
\end{Claim}

\begin{proof}
Throughout this phase, we will construct an increasing sequence of sets $\emptyset = A_{t_2}\subset \ldots \subset A_{t_3}$ with a triple separation. As during the second phase, we will show that if at some point $t$ it is Dominator's turn and the phase is not over, at least one of the following holds.

\begin{itemize}
    \item D can play the next move so that $\beta(t+2)\geq \beta(t)+3$, $\delta(t+2)\geq \delta(t)+1$ and $\nu(t+2)\leq\nu(t)$.
    \item D can play the next two moves so that $\beta(t+4)-\beta(t)\geq 6$, $\delta(t+4)\geq \delta(t)+1$ and $\nu(t+4)\leq\nu(t)$.
    \item D can play the next move so that $\beta(t+2)\geq \beta(t)+3$, $\chi(t+2)\geq \chi(t)+2$ and $\nu(t+2)\leq\nu(t)+1$.
    \item D can play the next move so that $\beta(t+2)\geq \beta(t)+3$, $\chi(t+2)\geq \chi(t)+2$, $\nu(t+2)\leq\nu(t)+2$ and $\vert A_{t+2}\vert \geq \vert A_t\vert +2$.
\end{itemize}

Assuming that this is true, it is easy to deduce the Claim inductively. In the following, all dominated white vertices are to be colored black immediately without explicit mention.

We now show that Dominator can always play as propositioned above. If there is a terminal walk, we instruct D to play as in the previous phase, leading to the first outcome.

Suppose next there was a white vertex $u$ with three neighbors $v_1, v_2$ and $v_3$ and that $v_1$ has another white neighbor, which we call $w_1$. Because no terminal walk exists, $v_2$ and $v_3$ must both also have a second white neighbor and we denote them by $w_2$ and $w_3$ respectively. Since no 2-circuit exists, $w_1,w_2$ and $w_3$ are distinct. As $G$ has minimum degree at least 2, $w_1$, $w_2$ and $w_3$ must all have a second neighbor, and because no vertex has at least three white neighbors, at least two of them must be distinct. Without loss of generality, we can therefore assume that $w_1$ and $w_2$ have distinct neighbors $x_1$ and $x_2$ respectively. Note that $x_1$ and $x_2$ cannot be adjacent to $u$. Going out from this setup we have to distinguish between a few different cases. For each case, we will only state the moves that D will play, and we leave it up to the reader to verify that we arrive at one of the four outcomes listed above.

Suppose first that $x_1$ and $x_2$ are both adjacent to $w_3$. In this scenario, D can play $v_1$. If S replies by playing a vertex that is not adjacent to $w_2$ or $w_3$, D can play $x_2$ in the next move, and we will be able to mark both $v_3$ and $v_2$ as depleted. Similarly, if S plays a vertex that is not $v_2$ but is adjacent to $w_2$ or $w_3$, we will be able to mark $v_2$ or $v_3$ as depleted immediately. This leaves the possibility that S replies by playing $v_2$. In this case, we mark $x_1,x_2$ and $v_3$ as dependent, noting that this increases $\sigma$ by three and $\nu$ only by one.

Next we assume that $x_1$ is adjacent to $w_3$ but $x_2$ is not. Since no terminal walk exists, $x_2$ must have a different white neighbor $y$. If $y=w_1$, the situation is isomorphic to that above as can be seen by relabeling $w_1$ as $w_3$ and vice versa. Assume therefore that $y\neq w_1$. Because $y$ is not a leaf and no vertex has three white neighbors, it has a second neighbor $z$ distinct from $v_1,v_2,v_3,x_1$ and $x_2$. Again, $z$ must have a second white neighbor $a$ different from $y$, and $a$ must be one of $w_1, w_3$ or $u$. Indeed, if $a$ was not, $azyx_2w_2v_2uv_3w_3x_1w_1v_1$ would be a 6-walk. 

Suppose that $a=u$ (see \Cref{fig:au}). In this case, D can play $x_2$. If S replies by playing a vertex that is adjacent to $u$, we can mark either $v_2$ or $z$ as depleted. If S plays $x_1$, we can mark $v_1,v_2,v_3$ and $z$ as dependent, thus increasing $\chi$ by at least three while increasing $\nu$ by at most one. If S does not play $x_1$ or a vertex that is adjacent to $u$, D can next play one of $v_1$ and $v_3$ and we mark both $v_2$ and $z$ as depleted. The cases $a=w_1$ and $a=w_3$ are symmetric so we will only consider $a=w_3$ (see \Cref{fig:aw}). D can still play $x_2$, and if S does not reply by playing a vertex that is adjacent to $u$ or $w_3$, D will play $v_3$ next and we can mark both $v_2$ and $z$ as depleted. The cases of S playing a vertex adjacent to $u$ or $w_3$ are again symmetric, so we may assume that S plays a vertex adjacent to $u$. If the vertex that S has played is not $v_2$, then we can mark $v_2$ as depleted, and if it is, D can play $x_1$ next upon which we mark all of $v_1,v_2$ and $z$ as depleted.

\begin{figure}
\centering
\begin{minipage}{.5\textwidth}
  \centering
  \includegraphics[width=.8\linewidth]{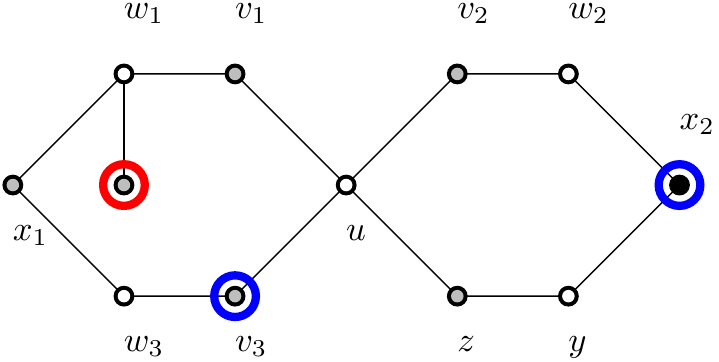}
  \captionof{figure}{The case where $a=u$.}
  \label{fig:au}
\end{minipage}
\begin{minipage}{.5\textwidth}
  \centering
  \includegraphics[width=.8\linewidth]{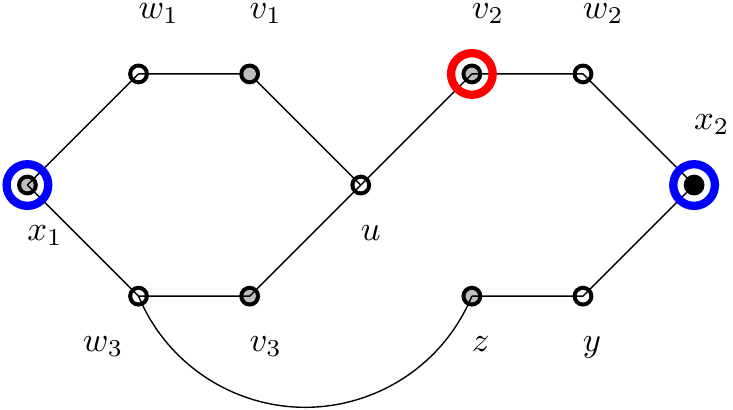}
  \captionof{figure}{The case where $a=w_3$.}
  \label{fig:aw}
\end{minipage}
\end{figure}

We are now left with the possibility that neither $x_1$ nor $x_2$ are adjacent to $w_3$ and have other neighbors $y_1$ and $y_2$ respectively. We will assume first that $y_1$ and $y_2$ are distinct. In this case, they must have neighbors $z_1,z_2\notin \{v_1,v_2,v_3,x_1,x_2\}$ respectively. Note that $z_1\neq z_2$ as otherwise $w_3v_3uv_1w_1x_1y_1zy_2x_2w_2v_2$ would be a 6-walk. In this situation, D plays $v_3$ and in the next move, depending on Staller's reply, one of $x_1$ and $x_2$ to mark $v_1$ or $v_2$ as depleted.

Assume that $y:=y_1=y_2$. Let $x_3\notin \{x_1,x_2\}$ be the second neighbor of $w_3$. Once more, $x_3$ must have a second white neighbor $a$, which must be different from $u$. On the other hand, we must have $a\in \{w_1,w_2,y\}$, as otherwise there would be a 6-walk. If $a=w_1$, we have a situation that is isomorphic to an earlier one as can be seen by relabeling $x_1$ as $z$, $x_3$ as $x_1$ as well as swapping the labels of $v_1$ and $w_1$ with those of $v_3$ and $w_3$ respectively (see \Cref{fig:aw}). An analogous isomorphism can be found if $a=w_2$. 

Finally, we have to consider the case where $a=y$. Here, D can play $v_1$. If S replies by playing a vertex that is not adjacent to $y$, D can play one of $x_2$ or $x_3$ next, marking $x_1$ as depleted in either case. If S replies by playing a vertex that is adjacent to $y$ but is not $x_1$, we can mark $x_1$ as depleted immediately. Lastly, if S replies by playing $x_1$, we can mark $v_2,v_3,x_2$ and $x_3$ as dependent, set $B_u=\{v_1,v_2,v_3\}$, $B_y=\{x_1,x_2,x_3\}$ and let $A_{t+2}=A_t\cup \{u,y\}$.

\begin{figure}[htbp]\centering
    			\includegraphics{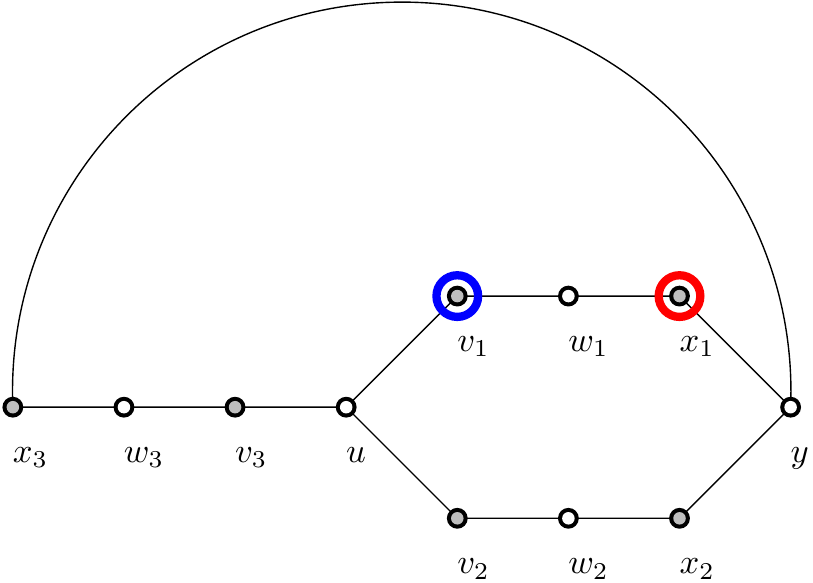}
    			\caption{The case where $y=a$.}
    			\label{fig:separated}
\end{figure}

If neither a terminal walk nor a vertex like $u$ exists, but the phase is not over, there must be a 5-walk $w_1\ldots v_5$. As $w_4v_4w_5v_5$ is not a terminal walk, $v_5$ must have a second white neighbor $w_6$ and in fact, we must have $w_6=w_1$. Indeed, because $w_2,w_3$ and $w_4$ have only two neighbors, $w_6$ must be distinct from them. If we also had $w_6\neq w_1$, $w_6$ would have another neighbor $v_6\notin \{v_1,\ldots v_5\}$ and we would have found a 6-walk, which cannot be. 

Hence, $w_1\ldots v_5$ is a 5-circuit and D can play $v_1$. If S replies by playing $v_2$, D plays $v_4$ next and we mark $v_3$ and $v_5$ as depleted. If S replies by playing $v_3$ or $v_4$, we can mark $v_2$ or $v_5$ respectively as depleted immediately. Lastly, if S plays neither of these vertices, $w_3$ and $w_4$ are still white after Staller's move and D can play $v_3$ after which we mark $v_2$ as depleted. 
\end{proof}

At the end of Phase 3, no terminal walk, no 5-walk and no white vertex with degree higher than two and part of a 2-walk exists, and thus the remaining white vertices must lie in very specific configurations, as the following lemma shows.

\begin{lemma}\label{lemma:circuit-3-4}
Once Phase 3 ends, every vertex with two white neighbors lies inside a 3- or 4-circuit.
\end{lemma}
\begin{proof}
Suppose a vertex $v_1$ has two white neighbors $w_1$ and $w_2$. Because $w_2$ is not a leaf, $w_2$ has an additional neighbor $v_2$ and since no terminal walk or 2-circuit exists, $v_2$ has another white neighbor $w_3$ which must in turn have a second neighbor $v_3\neq v_1$. Repeating this argument, we know that $v_3$ has a second white neighbor $w_4$ that has itself a second neighbor $v_4\neq v_2$. If $w_4=w_1$ (and consequently $v_4=v_1$), $v_1$ is in a 3-circuit, so we will assume that $w_4\neq w_1$. Repeating the argument again, $v_4$ has a second white neighbor $w_5$ with a second neighbor $v_5$. We know that $w_5\neq w_2$ since $w_2$ can only have two neighbors, and therefore we must have $w_5=w_1$ as otherwise we would have found a 5-walk. Thus, $v_1$ is in a 4-circuit.
\end{proof}

At this point, it is possible that there are vertices $v,w,x$ such that $w$ is undominated and adjacent to both $v$ and $x$, and $v$ is marked as dependent but $x$ is not. In these situations, $x$ cannot have a second white neighbor, as that would yield a terminal walk, and $x$ cannot be marked as depleted by definition. Furthermore, by marking all such $x$ that are neighbors of white neighbors $w$ of dependent vertices $v$ as dependent, we do not increase $\nu$. After doing so, we can make the following observation.

\begin{lemma}\label{lemma:nu-chi}
We have $\nu(t_3)\leq \chi(t_3)$.
\end{lemma}

\subsection{Phase 4}
Phase 4 ends once the game is over.

\begin{Claim}\label{phase:endgame}
Let $C(t)$ be the set of vertices that have a single white neighbor after move $t$. Dominator can play in a way so that there exists non-negative integers $a$ and $b$ such that the following hold.

\begin{itemize}
    \item $T_4\leq 2a+b+1$.
    \item $\beta(t_4)-\beta(t_3) \geq 3a+b$.
    \item $\delta(t_4)-\delta(t_3) \geq a+b-\nu(t_3)$.
    \item For all $t\geq t_3$, $C(t)\subset C(t_3)$
    \item No new vertices are marked as dependent.
    \item Every vertex with two white neighbors lies inside a 3- or 4-circuit.
\end{itemize}
\end{Claim}

Note that the last point holds at the start of the phase by \Cref{lemma:circuit-3-4}.

\begin{proof}
We denote by $b$ the number of moves in which a vertex in $C(t_3)$ is played. Whenever either player does so, we paint its unique white neighbor $w$ black. Because $w$ is not a leaf, it has a second neighbor $x$. If $x$ is marked as dependent, $\nu$ has certainly decreased by one. If $x$ is not marked as dependent, we can mark it as depleted. Overall, the number of depleted vertices introduced from such moves must be at least $b-\nu(t_3)$.

We can now give Dominator's strategy for this phase. If $C(t_3)\neq \emptyset$, then on the first move there has to exist a vertex $v$ with a single white neighbor $w$. After D plays $v$, we can paint $w$ black. Note also that $C(t_3+1)\subset C(t_3)$. If $C(t_3)=\emptyset$, but the game is still ongoing, D plays an arbitrary vertex with two white neighbors and we do not paint or mark any vertices.

From here on out Dominator plays reactively depending on Staller's previous move. We distinguish between the following cases.

\begin{itemize}
    \item S plays a vertex $v_i$ inside a 4-circuit $w_1\ldots v_4$, and without loss of generality, $i=1$. After D responds by playing $v_3$, we paint all of $w_1,\ldots, w_4$ black and mark $v_2$ and $v_4$ as depleted. Because all of $w_1,\ldots,w_4$ have degree 2, $C(t+2)=C(t)$. We have $\beta(t+2)=\beta(t)+4$ and $\delta(t+2)=\delta(t)+2$.
    \item S plays a vertex $v_i$ inside a 3-circuit $w_1\ldots v_3$, and without loss of generality, $i=1$. After D responds by playing $v_2$, we paint $w_1,w_2$ and $w_3$ black and mark $v_2$ and $v_4$ as depleted. Again, $C(t+3)=C(t)$. We have $\beta(t+2)=\beta(t)+3$ and $\delta(t+2)=\delta(t)+1$.
    \item S plays a vertex in $v\in C(t)$. If $C(t+1)$ contains another vertex $u$, D plays it and we paint its white neighbor black as well. If $C(t+1)$ is empty, but the game is not over, D plays an arbitrary vertex with two white neighbors and we do not paint or mark any vertex. Note that this situation can occur at most once and only if $C(t_3)\neq \emptyset$.
\end{itemize}

Denoting the combined number of moves played by Staller of the first type and second type by $a$, we obtain the Claim.
\end{proof}

\subsection{Bounding the number of moves}

Let $T=t_4=T_1+T_2+T_3+T_4$. To prove \Cref{Bound_Min_degree_2}, it is enough to show that $T\leq 5(n+1)/7$. We will do so by combining the bounds on $\beta,\delta$ and $\chi$ that we have established during the previous sections. 
As mentioned before, we have two principle angles of attack. Firstly, we know that the number $\beta(T)$ of black vertices at the end of the game is at most $n$. Summing the lower bounds on the number of vertices that have been painted black during each of the phases, we get that
\begin{align}\label{eq:blackdeg2}
     2T_1+\frac{3T_2}{2}+\frac{3T_3}{2}+3a+b\leq \beta(T)\leq n
\end{align}

Secondly, we can divide the vertices up in those that have been played and those that have not been played. Of the former kind there are exactly $T$ many. The set of not-played vertices will include all $\delta(T)$ depleted vertices and at least $\chi(T)$ dependent vertices. Therefore, we have

\begin{align}\label{eq:min-deg2-played-unplayed}
    T+\delta(t_4)+\chi(t_4)\leq n.
\end{align}
By the bounds on $\delta$ and $\chi$ we observed in \Cref{phase:long-walks}, \Cref{phase:intersecting-circuits}, \Cref{phase:endgame} and the monotonicity of $\chi$, we have
\begin{align}\label{eq:LBdelta-chi}
    \delta(t_4)+\chi(t_4)\geq a+b-\nu(t_3)+\frac{T_3}{4}+\frac{3\chi(t_3)}{4}+\frac{T_2}{4}.
\end{align}

Adding the bounds $\nu(t_3)\leq \chi(t_3)$ from \Cref{lemma:nu-chi} and $\nu(t_3)\leq (\chi(t_3)+M)/2$ from \Cref{phase:intersecting-circuits} we get $\nu(t_3) \leq 3\chi(t_3)/4+M/4$. Inserting this bound and $T_1\geq M$ from \Cref{phase:3+} into \ref{eq:LBdelta-chi}, we get
\begin{align}\label{eq:min-degree-2-delta-chi}
    \delta(t_4)+\chi(t_4)\geq a+b+\frac{T_2+T_3}{4}-\frac{T_1}{4}.
\end{align}

Adding \eqref{eq:blackdeg2} to $2/3$ times \eqref{eq:min-deg2-played-unplayed}, we obtain
\begin{align*}
    \frac{8T_1}{3}+\frac{13(T_2+T_3)}{6}+\frac{2(T_4+\delta(t_4)+\chi(t_4))}{3}+3a+b\leq \frac{5n}{3}.
\end{align*}
Combining this with \eqref{eq:min-degree-2-delta-chi} we get
\begin{align*}
    \frac{5T_1}{2}+\frac{7(T_2+T_3)}{3}+\frac{2T_4+11a+5b}{3}\leq \frac{5n}{3}.
\end{align*}
But, by \Cref{phase:endgame}, we have $T_4\leq 2a+b+1$, which gives $7T/3-5/3\leq 5n/3$ and thus $T\leq 5(n+1)/7$. This proves \Cref{Bound_Min_degree_2}.

\section{Proof of \Cref{3-4_thm}}

Suppose towards a contradiction that \Cref{3-4_thm} were false. For the entirety of this section, we fix a counterexample $G$ that is a minimal with respect to the number of its edges. As before, let $n$ be the number of vertices of $G$.

\subsection{Simplifications}

\begin{lemma}\label{lemma:no-two-leaves}
No parent in $G$ has more than one leaf.
\end{lemma}
\begin{proof}
Suppose $G$ has an owner $x$ with two leaves $y$ and $z$. By minimality of $n$, D must have a strategy $\Delta$ for the game on $G-z$ to end in $3n/4$ moves. Since $y$ and $z$ can never both be played in the game on $G$, D can simply follow $\Delta$ for the game on $G$, pretending when S plays the vertex $z$ that they played $y$ instead.
\end{proof}

\begin{lemma}\label{jans-lemma}
Let $u$ be a parent in $G$ with a leaf $w$ and another neighbor $v$. Then $G-uv$ has an isolated edge.
\end{lemma}
\begin{proof}
We know that $G-uv$ cannot have an isolated vertex as $v$ cannot be another leaf of $u$, and we assume for contradiction that $G-uv$ has no $K_2$ component either. By minimality of $G$, we know that D has a strategy $\Pi$ to play for a dominating set with at most $3n/4$ vertices in $G-uv$. We claim that this must be true in $G$ as well.

Indeed, the only way D could fail to apply $\Pi$ in $G$ is that S plays a vertex that would be forbidden in $G-uv$ because all its neighbors are already totally dominated. The only vertices that have additional neighbors in $G$ are $u$ and $v$, but $u$ can always be played because of its leaf. If on the other hand S can play $v$ in $G$ but not in $G-uv$, we must have that firstly $w$ has not yet been played and secondly all neighbors of  $v$ except $u$ have already been otherwise totally dominated. But then D may pretend that S plays $w$ instead of $v$ as this is still allowed in $G-uv$, and by doing so D can play according to $\Pi$ for the rest of the game.
\end{proof}

In the rest of this section, we will make the following distinction between leaves. A \emph{leaf of type A} is a leaf whose parent is not adjacent to another parent, and a \emph{leaf of type B} is a leaf whose parent is adjacent to at least on other parent. \\

Let $x$ be a leaf of type B, and suppose that the parent $u$ of $x$ is adjacent to $k\geq 1$ other parents $v_1,\ldots,v_k$ with leaves $y_1,\ldots,y_k$ respectively. Then we remark that by \Cref{jans-lemma}, the component of $G$ containing $u$ consists only of $C=G[\{u,x,v_1,y_1,\ldots,v_k,y_k\}]$. \\

On the other hand, if $x$ is a leaf of type A, then the parent $y$ of $x$ has exactly one other neighbor $z$. By \Cref{lemma:no-two-leaves}, $z$ is not a leaf and by definition of a leaf of type A, $z$ is not a parent. Abusing language, we say that $z$ is the \emph{grandparent} of both $x$ and $y$. 
Let $L$ denote the set of leaves (of either type), let $P$ denote the set of parents and let $F$ denote the set of grandparents. Furthermore, let $R = V(G) \setminus (L \cup P \cup F)$. \\

As in the minimum degree at least 2 case, the strategy is divided into consecutive \emph{phases}, each of them starting after the last one ends by meeting a predefined end condition. We track the number of moves played so far by a variable $t$. There are going to be six phases, and for each $i\in [6]$, $t_i$ is the last move before the end of Phase $i$ and let $t_0=0$. It can happen that phases are skipped because their end conditions are met from the start of the phase and in this case we have $t_i=t_{i-1}$. For each $i\in [6]$, let $T_i=t_i-t_{i-1}$.

\subsection{Phase 1}

The first phase starts at the start of the game and ends once it is Dominator's turn, no grandparent is adjacent to two or more white parents and no parent is adjacent to one or more white parent.

\begin{Claim}\label{AnalysisPhase1}
D can play in a way such that there exists some non-negative integers $r_F$ and $r_P$ such that the following hold.
\begin{itemize}
    \item $\beta(t_1)\geq \frac{3}{2}T_1$.
    \item $\delta(t_1)\geq 2r_F+r_P$.
    \item $r_F+r_P \geq T_1/2$.
    \item $\lambda(t_1)\leq \delta(t_1)+T_1/2$.
    \item No vertex has been marked as dependent.
    \item Each leaf of type B has been played or marked as depleted.
\end{itemize}
\end{Claim}

\begin{proof}
For as long as possible, D plays grandparents $v$ with two white neighbors in $P$, and we paint all white parents of $v$ in $P$ black and mark their leaves as depleted. If no such vertex $v$ exists, then D plays parents $w$ with at least one white neighbour $u$ in $P$, and we paint all white neighbours of $w$ in black and mark the leaf of $u$ as depleted.

With regard to moves by Staller, whichever vertex $z$ they play, we paint all neighbors of $z$ black. As $z$ must have had at least one undominated neighbor before Staller's move, we paint at least three vertices black every two moves and the first part follows. 

Let $r_F$ and $r_P$ be respectively the number of times that D plays a grandparent and a parent during Phase 1. The second part of the claim holds because we mark at least 2 or 1 leaves as depleted whenever D plays a grandparent or parent, respectively. As D makes at least $T_1/2$ moves, the third part follows immediately. The fourth point follows from the fact that D never plays a leaf during this phase and S plays at most $T_1/2$ moves.

Once the phase is over, we mark all leaves whose parents are painted black as depleted, from which the last point follows. Indeed, if a leaf $y$ of type B has neither been played nor marked as depleted, then its parent $u$ is white. As $y$ is a leaf of type B, $u$ must have a neighbor in $P$, which we call $v$, and since $u$ is white, we know that $v$ has not been played. Hence $v$ has a white neighbor in $P$, which contradicts that Phase 1 has ended.
\end{proof}

\subsection{Phase 2}

Phase 2 ends once no grandparent is adjacent to a white parent and at least three white vertices overall and it is Dominator's turn.

\begin{Claim}\label{AnalysisPhase2}
D can play in a way such that the following hold.

\begin{itemize}
    \item No vertex has been marked as dependent.
    \item $\lambda(t_2)-\lambda(t_1)\leq T_2$
    \item $\beta(t_2)-\beta(t_1) \geq 2T_2$
    \item $\delta(t_2)-\delta(t_1) \geq T_2/2$
\end{itemize}
\end{Claim}

\begin{proof}
If the phase is ongoing and D has to make a move, they play a grandparent $u$ that is adjacent to a white parent and at least three white vertices overall. After they have done so, we paint all white neighbors of $v$ black. One of these neighbors is a parent and we mark its leaf as depleted.

We handle Staller's moves as in the first phase and leave it to the reader to verify the Claimed bounds.
\end{proof}

\subsection{Phase 3}
During this phase, we will decide after each move of Staller whether or not to set a \emph{reaction flag}. The phase ends when it is Dominator's turn, the reaction flag is not set, there exists no parent so that both its leaf and grandparent are white and there exists no grandparent that is adjacent to a white parent and a second white vertex.

\begin{Claim}\label{AnalysisPhase3}
D can play in a way such that there exists some non-negative integers $m_1,m_2,m_3,p$ and $q$ such that the following hold.

\begin{itemize}
    \item No leaf has been marked as dependent.
    \item $\beta(t_3)-\beta(t_2) \geq \frac{3}{2}(T_3-T_2)$
    \item $\delta(t_3)-\delta(t_2) = m_1+m_3+q$
    \item $\chi(t_3) \geq m_2$
    \item $\lambda(t_3)-\lambda(t_2) = m_1+m_2+m_3+q$
    \item $T_3/2=m_3+p+q$
\end{itemize}
\end{Claim}

\begin{proof}
Throughout this phase, before each move by Staller we paint all white totally dominated vertices black. In particular, since a leaf can only be totally dominated if its parent has been played, a grandparent of a black leaf must be black.

If it is Dominator's turn and the reaction flag is not set, D will play a parent with a white leaf and a white grandparent where possible after which we color both the leaf and the grandparent black. If no such parent exists and the phase is still ongoing, there must be a grandparent $u$ that is adjacent to a white parent $x$ and some other white neighbor. Then D plays $u$ and we color the parent as well as the other neighbor black and mark the leaf of $x$ as depleted. For those cases where the reaction flag is set, we instruct D to play dynamically depending on Staller's previous move.

If S plays any vertex that is not a leaf, we do not set the reaction flag and simply paint all its white neighbors black. Whenever S does plays a leaf with a parent $x$ and a grandparent $y$, we categorize this move into one of three different types as follows. 

\begin{itemize}
    \item $y$ does not have any white neighbor other than $x$.
    \item $y$ has one other white neighbor $z\notin P$, and $z$ has no other neighbor having another white neighbor.
    \item $y$ has one white neighbor $z$, and $z$ has at least one neighbor $u$ with another white neighbor $v$.
\end{itemize}

Due to the end conditions of the previous phases, no grandparent is adjacent to more than one white parent, so the above list is exhaustive. Let $m_1,m_2,m_3$ be the respective number of moves of the above listed types played by S during Phase 3. Furthermore, we denote by $p$ the number of times that the reaction flag was not set and D played a parent and by $q$ the number of times the reaction flag was not set and D played a grandparent.

If Staller plays a leaf, we always color its parent black and proceed depending on the type $m_1, m_2$ or $m_3$ of the leaf. If the leaf has type $m_1$, we can mark its grandparent $y$ as depleted (note that $y$ cannot have been played before, as that would make it illegal for S to play the leaf) and do not set the reaction flag. 

Similarly, if S plays a leaf of type $m_2$, we consider the unique white neighbor $z\notin P\cup L$ of its grandparent $y$. Since $z$ is not a leaf, it must have some other neighbor $u$, and by definition of type $m_2$, $u$ has no other white neighbor than $z$. We can therefore mark $y$ as dependent and $u$ as well if it has not been marked as such already. In either case, we increase $\chi$ by at least 1. We do not set the reaction flag.

If S plays a leaf of type $m_3$, we do set the reaction flag and instruct D to play his next move as follows. For Staller's move to have type $m_3$ the grandparent $y$ has a white neighbor $z\notin P\cup L$ which in turn has a neighbor $u$ with another white neighbor $v$. Here, D plays $u$ and we can mark $z$ as depleted.
\end{proof}

Before moving on to the next phase, we make an important observation about the leaves of the graph. Let $I$ denote the set of leaves that have not been played or marked as depleted. As documented in \Cref{AnalysisPhase1}, all such leaves must be of type A, and by definition of $\lambda(t)$, $I$ has at least $\vert L\vert-\lambda(T_3)\geq |L|-m_1-m_2-m_3-q-T_2-\lambda(t_1)$ elements. Let now $J$ be the set of grandparents of leaves in $I$ and let $C\subset J$ be those among them that have not been played or marked as depleted or dependent. We now mark all vertices in $I\cup C$ as dependent. This increases $\chi$ by at least $|I|$. Indeed, if the grandparent of a leaf $x\in I$ is in $C$, they share the parent of $x$ as their only undominated neighbor. If the grandparent of $x$ is not in $I$, the parent of $x$ is either dominated or the grandparent of $x$ had already been marked as dependent at some earlier point. Therefore, we have the following.

\begin{lemma}\label{lemma:dependent-leaves}
We have $\chi(t_3)\geq |L|-m_1-m_3-q-T_2-\lambda(t_1)$.
\end{lemma}

\subsection{Phase 4}
Phase 4 ends once it is Dominator's turn and no vertex has three or more white neighbors.

\begin{Claim}\label{AnalysisPhase5}
D can play in a way such that the following hold.

\begin{itemize}
    \item $\beta(t_4)-\beta(t_3) \geq 2T_4$.
\end{itemize}

\end{Claim}

\begin{proof}
Dominator keeps playing vertices $y$ with at least three white neighbors, and after Staller plays a vertex $z$ with at least one white neighbor, we paint all neighbors of $z$ and $y$ black. The claimed bound follows.
\end{proof}

\subsection{Phase 5}

Phase 5 ends once it is Dominator's turn and no vertex has two white neighbors.

\begin{lemma}\label{AnalysisPhase6}
D can play in a way such that the following hold.

\begin{itemize}
    \item No vertex that has been marked as dependent ceases to be marked as dependent.
    \item $\beta(t_5)-\beta(t_4)\geq \frac{3}{2}T_5$
    \item $4(\delta(t_5)-\delta(t_4))+2(\chi(t_5)-\chi(t_4))\geq T_5$
\end{itemize}

\end{lemma}

\begin{proof}
Similar to Phases 2 and 3 of the proof of \Cref{Bound_Min_degree_2}, we will show that whenever it is Dominator's turn, at least one of the following hold.
\begin{itemize}
    \item D can play the next move so that
    \begin{align*}
        \beta(t+2)-\beta(t)\geq 3 \qquad \text{and} \qquad \delta(t+2)-\delta(t)+\chi(t+2)-\chi(t)\geq 1,
    \end{align*}
    \item D can play the next two moves so that 
    \begin{align*}
    \beta(t+4)-\beta(t)\geq 6 \qquad \text{and}\qquad
    \delta(t+4)-\delta(t)\geq 1.
    \end{align*}
\end{itemize}
Once we have established this dichotomy, the claim follows again from a simple inductive argument.

Note that while $G$ can now have leaves, no leaf can be part of any $k$-walk for $k\geq 2$ at this point. Indeed, if we take any $k$-walk $w_1\ldots v_k$, only $w_1$ and $v_k$ can possibly be leaves. But if $w_1$ were a leaf, then $v_1$ would be a parent with a white leaf and a white neighbor, which cannot exist after the end of Phase 3. Likewise, if $v_k$ were a leaf, $v_{k-1}$ would either be a parent or a grandparent with a white parent and a second white neighbor, which cannot exist after the end of Phase 3 either. 

For this reason, D can play first as in \Cref{sec:phase2-mindeg2} and then as in \Cref{sec:phase3-mindeg2} until, as by \Cref{lemma:circuit-3-4} in \Cref{sec:phase3-mindeg2}, every vertex with two white neighbors lies in a 3- or 4-circuit.

Suppose now that it is Dominator's turn and there is 3-circuit $w_1\ldots v_3$. Here, D can play $v_1$, we paint $w_1$ and $w_2$ black and mark $v_2$ and $v_3$ as dependent. Should there be instead a 4-circuit $w_1\ldots v_k$, D can still play $v_1$, but we need to consider two different continuations. If S plays a vertex that is adjacent to $w_3$ or $w_4$ next, without loss of generality $w_3$, we paint $w_1,w_2$ and $w_3$ black and mark $v_3$ and $v_4$ as dependent. If S plays a vertex that is not adjacent to either of them, D can play $v_3$ next, after which we paint $w_1,w_2,w_3$ and $w_4$ black and mark $v_2$ and $v_4$ as depleted.
\end{proof}

\subsection{Phase 6}
Phase 6 ends once no undominated vertex is left. From this point on, we paint all dominated vertices black immediately. Also, we mark all vertices that have one white neighbor as dependent.

\begin{Claim}\label{AnalysisPhase7}
Regardless of what moves either player makes, we have $\beta(t_6)-\beta(t_5) \geq T_6$ and $T_6\leq \chi(t_5)+L-p-r_P$.
\end{Claim}

\begin{proof}
Every vertex that is played during this phase is either dependent or dominates a leaf. There are at most $L-p-r_P$ undominated leaves left at the start of the phase and since each neighbor of dependent parents can only be dominated once and no new dependent vertices are introduced, at most $\chi(t_5)$ dependent vertices can be played.
\end{proof}

\subsection{Bounding the number of moves}

Let $T=T_1+\ldots+T_6$. To prove \Cref{3-4_thm}, it is enough to show that $T\leq 3n/4$. We will do so by combining the bounds on $\beta,\delta,\lambda$ and $\chi$ we have established during the previous sections.

As in Section 3, we sum the lower bounds on the number of new black vertices from each of the phases to obtain
\begin{align}\label{eq:black}
     3(T_1+T_3+T_5)/2 + 2(T_2+T_4) + T_6 \leq \beta(T)\leq n.
\end{align}
We also divide the vertices up in those that have been played and those that have not been played, which gives
\begin{align}\label{eq:played-unplayed}
    T+\delta(T)+\chi(T) \leq n.
\end{align}
Moreover, we combine our bounds for the increases of $\delta$ and $\chi$ during the different phases. For $\delta$, we obtain
\begin{align}\label{eq:u1}
    \delta(T)\geq \lambda(t_1)-T_1/2+T_2/2+m_1+m_3+q+\delta(t_5)-\delta(t_4)
\end{align}
or by using that $m_3+p+q=T_3/2$ from \Cref{AnalysisPhase3} and $\delta(t_1) \geq 2r_F+r_P$ from \Cref{AnalysisPhase1},
\begin{align}\label{eq:u2}
    \delta(T)&\geq 2r_F+r_P+T_2/2+T_3/2-p+\delta(t_5)-\delta(t_4).
\end{align}
Using \Cref{lemma:dependent-leaves}, we have
\begin{align}\label{eq:x1}
    \chi(T)\geq |L|-m_1-m_3-q-T_2-\lambda(t_1)+\chi(t_5)-\chi(t_4).
\end{align}
Recall also that we have from \Cref{AnalysisPhase7} the inequality
\begin{align}\label{eq:x2}
    \chi(T)\geq T_6-|L|+p+r_P.
\end{align}
Adding up the inequalities \eqref{eq:u1}, \eqref{eq:u2}, \eqref{eq:x1} and \eqref{eq:x2} and using that $2(r_P+r_F) \geq T_1$, we obtain
\begin{align*}
    2\delta(T)+2\chi(T)\geq T_1/2+T_3/2+T_6+2(\delta(t_5)-\delta(t_4))+\chi(t_5)-\chi(t_4).
\end{align*}
By \Cref{AnalysisPhase6}, $2(\delta(t_5)-\delta(t_4))+\chi(t_5)-\chi(t_4)\geq T_5/2$, so we must have
\begin{align}\label{eq:u-x}
    2\delta(T)+2\chi(T)\geq T_1/2+T_3/2+T_5/2+T_6.
\end{align}
To get a bound on $T$, we linearly combine four times \eqref{eq:played-unplayed} together with two times \eqref{eq:black} to obtain
\begin{align*}
    7T_1+8T_2+7T_3+8T_4+7T_5+6T_6+4\delta(T)+4\chi(T)\leq 6n,
\end{align*}
and by inserting \eqref{eq:u-x} into this we arrive at $8T\leq 6n$, which completes the proof.\qed

\section{Concluding remarks and open problems}

At the first glance, it is perhaps surprising that the bounds in \Cref{3-4_thm} and \Cref{Bound_Min_degree_2} are so close to clean fractions of $n$, rather than having larger error-terms. However, it turns out that this is part of a pattern, dictating that for a large class of variants of the domination game, the error term of upper bounds as in \Cref{3-4_thm} and \Cref{Bound_Min_degree_2} can be at most 1.

To make this precise, we consider a variant of the total domination game, that can be seen as a generalization of most other variants that have been studied. The \emph{transversal game} was defined in \cite{bujtas2016transversal} as a game played on a hypergraph $H$ where the two players, Edge-hitter and Staller, alternately select a vertex from $H$, with the rule that each newly selected vertex must \emph{hit}, i.e., be contained in, at least one edge that does not intersect the set of previously selected vertices. The game ends when the set of selected vertices becomes a transversal in $H$, i.e., every edge of $H$ intersects the set of selected vertices. Edge-hitter aims to end the game in as few moves as possible, whereas Staller wants the game to last as long as possible. The \emph{game transversal number} (respectively \emph{Staller-start game transversal number}) $\tau_g(H)$ (respectively $\tau'_g(H)$) of $H$ is the number of moves played if Edge-hitter (respectively Staller) starts the game and both players play optimally.

\begin{lemma}\label{copylemmageneral}
Let $\cC$ be a set of hypergraphs closed under taking multiple disjoint copies of a hypergraph in $\cC$ and suppose that there exists $c>0$ such that for every hypergraph $H \in \cC$ on $n$ vertices, $\tau_{g}(H) \leq (c+o(1))n$. Then we also have for every $H\in \cC$ on $n$ vertices, $\tau_{g}(H) \leq cn+1$. 
\end{lemma}

\begin{proof}
For a given hypergraph $H$ and a set of vertices $A$ in $H$, let $T(H;A)$ be the \emph{transversal game in $H$ starting from $A$}, in which the players play the Edge-Hitter-start transversal game starting with the set $A$ already selected. Likewise, let $T'(H;A)$ be the analogous game in which Staller starts. We define $\tau_g(H;A)$ and $\tau_g'(H;A)$ to be the number of selected vertices excluding $A$ at the end of the respective games $T(H;A)$ and $T'(H;A)$ if both players play optimally.

Note that for any $H$ and $A$, we have $\tau_g(H;A)\leq \tau_g'(H;A)+1$. Indeed, suppose that D has the legal move $v$ in the game $T(H;A)$. If $\tau_g'(H;A\cup \{v\})=0$, then $\tau_g(H;A)=1$ and there is nothing left to show. Otherwise, there must be a vertex $w$ such that $\tau_g(H;A\cup \{v,w\})\geq \tau_g(H;A)-2$. But by a straightforward generalization of the \emph{continuation principle} (see \cite{kinnersley2013extremal}), we have $\tau_g(H;A\cup \{v,w\})\leq \tau_g(H;A\cup\{w\})$, which is in turn at most $\tau_g'(H;A)-1$.

In order to arrive at a contradiction, assume that there exist $e>1$ and $H\in \cC$ on $n$ vertices such that $\tau_g(H)> cn+e$. Let $H^m$ be the graph consisting of $m$ disjoint copies $H_1,\ldots,H_m$ of $H$. For a set $A\subset V(H^m)$, let $k(A)$ be the number of copies $H_i$, in which there is a playable vertex in the game $T(H_i;A\cap V(H_i))$. Note that $k(A)$ is decreasing in $A$. Abbreviating $\tau_g(H_i;A\cap V(H_i))$ to $\tau_g(H_i;A)$, we say that a set $A$ is \emph{good}, if
\begin{align*}
    \tau_g(H^m;A)\geq -k(A)+\sum_{j\in [m]} \tau_g(H_i;A).
\end{align*}
If $\emptyset$ were good, we would have
\begin{align*}
    \tau_g(H^m)=\tau_g(H^m;\emptyset)\geq m\tau_g(H;\emptyset)-m>cmn+(e-1)m.
\end{align*}
Since $H^m$ has $mn$ vertices, this contradicts our assumption if $m$ is large enough. We may therefore assume that for some sufficiently large fixed $m$, $\emptyset$ is bad, and we consider a bad set $A$ of maximal size in $V(H^m)$. Clearly, $\tau_g(H^m;A)\geq 2$ and there must be a vertex $v$ such that for $A'=A\cup \{v\}$, we have $\tau_g'(H^m;A')= \tau_g(H^m;A)-1$. Let $i\in [m]$ be such that $v\in V(H_i)$. If $\tau_g'(H_i;A')>0$, then there must be a vertex $w\in V(H_i)$ such that for $A''=A'\cup\{w\}$, $\tau_g(H_i,A'')=\tau_g(H_i,A)-2$. Since $A''$ is good and $A$ is bad, we have
\begin{align*}
    \tau_g(H^m;A'')&\geq -k(A'')+\tau_g(H_i,A'')+\sum_{\substack{j\in [m]\\j\neq i}} \tau_g(H_j;A)\\
    &\geq -k(A)-2+\sum_{j\in [m]} \tau_g(H_j;A)\\
    &>\tau_g(H^m;A)-2,
\end{align*}
but as $\tau_g(H^m;A'') \leq \tau'_g(H^m;A')-1$, this contradicts $\tau_g'(H^m;A')= \tau_g(H^m;A)-1$.

If we have instead that $\tau_g'(H_i;A')=0$, then $k(A)=k(A')+1$. Let $r\in [m]$ be such that $\tau_g'(H_r,A')>0$. We know that $\tau_g'(H_r,A')\geq \tau_g(H_r,A')-1$, and thus there must be a $w\in V(H_r)$ such that for $A''=A'\cup\{w\}$, we have $\tau_g(H_r,A'')\geq\tau_g(H_r,A')-2$. But since $A$ is bad and $A''$ is good, we have
\begin{align*}
    \tau_g(H^m;A)&< -k(A)+\sum_{j\in [m]} \tau_g(H_j;A)\\
    &= -k(A')+\sum_{j\in [m]} \tau_g(H_j;A')\\
    &\leq -k(A'')+2+\sum_{j\in [m]} \tau_g(H_j;A'')\\
    &\leq \tau_g(H^m;A'')+2,
\end{align*}
contradicting again $\tau_g'(H^m;A')= \tau_g(H^m;A)-1$.
\end{proof}

We remark that the remaining error term of 1 can be removed, if the bound $\tau_g(H)\leq (c+o(1))n$ can be achieved in a version of the transversal game where Staller is allowed to forfeit his move whenever they wish. To see how \Cref{copylemmageneral} relates to other variants of the domination game, we take the total domination game as an example. The open neighborhood hypergraph of a graph $G$ is defined as the hypergraph $H_G$ with vertex set $V(H_G) = V(G)$ and hyperedge set $E(H_G) = \{ N_G(x) | x \in V(G) \}$ consisting of the open neighborhoods of vertices in $G$. As remarked in \cite{bujtas2016transversal}, the total domination game played on a graph $G$ can be seen as a special instance of the transversal game being played on the open neighborhood hypergraph $H_G$ of $G$. Since the class $\{H_G:G \text{ is a graph}\}$ is closed under taking disjoint copies, \Cref{copylemmageneral} also applies to the total domination game. \\

As mentioned earlier, \Cref{3-4_thm} is sharp, as the example of disjoint copies of $P_4$ or $P_8$ shows. However, it is unclear whether \Cref{3-4_thm} is asymptotically sharp if we add the condition that $G$ is a connected graph, and Bre{\v{s}}ar, Henning and Klav{\v{z}}ar in \cite{brevsar2021domination} raised the following problem.

\begin{Problem}
    Determine all connected graphs satisfying $\gamma_{tg}(G) = \frac{3}{4}n$.
\end{Problem}

So far the largest graph $G$ known to satisfy $\gamma_{tg}(G) = \frac{3}{4}n$ is a tree on 12 vertices constructed in \cite{brevsar2021domination}. \\

With respect to \Cref{Bound_Min_degree_2}, we first remark that the strategy described in this paper cannot give a better bound than $\frac{5}{7}n$. Indeed, if $G$ is the union of $k$ copies of $C_7$ and D follows the strategy laid out in Section 3, then it is easy to see that S has a strategy so that the game will last exactly $5k$ moves. However, even though $\gamma_{tg}(C_7)=5$, for the disjoint union $H_k$ of $k$ copies of $C_7$, we have that $\gamma_{tg}(H_k)=4k+1$, which provokes the question whether \Cref{Bound_Min_degree_2} is asymptotically tight. More generally, we ask the following question, which is analogous to the $1/2$-Conjecture made in \cite{bujtas20221}.

\begin{Question}\label{questionmindeg2}
What is the smallest constant $c>0$ such that for any graph $G$ with minimum degree at least $2$ we have $\gamma_{tg}(G) \leq cn+1$?
\end{Question}

Note that if $G$ is a disjoint union of triangles, then $\gamma_{tg}(G)=2n/3$. Thus, $c$ must be at least $\frac{2}{3}$, and it seems plausible that $c=\frac{2}{3}$ is in fact the answer to \Cref{questionmindeg2}. We remark further that by \Cref{copylemmageneral}, it is sufficient to show the bound $\gamma_{tg}(G)\leq (c+o(1))n$.\\

The domination game itself was introduced by Bre{\v{s}}ar, Klav{\v{z}}ar and Rall in \cite{brevsar2010domination}. The difference between the standard and the total domination game is that a vertex selected by the players \emph{dominates} itself, whereas is does not \emph{totally} dominates itself. The \emph{the game domination number} $\gamma_{g}(G)$ of $G$ is defined as number of vertices selected when Dominator starts the game and both players play optimally. One of the central topics related to the domination game is the 3/5-conjecture posed by Kinnersley, West and Zamani \cite{kinnersley2013extremal}.

\begin{conjecture}\label{3-5_conj}
If $G$ is an isolate-free graph of order $n$, then $\gamma_g(G) \leq \frac{3}{5}n$.
\end{conjecture}

A considerable amount of research has been done on the domination game and \Cref{3-5_conj} in particular \cite{kinnersley2013extremal,bujtas2015domination,schmidt20163,marcus2016domination,henning2016domination,BujtasMinDegree3,bujtas2020general,bujtas20221,versteegen2022domination}.

There are also many other variants of the domination game, all of which are studied at the moment, and we want to bring some of them to the attention of the reader. For instance, there are the conncected \cite{borowiecki2019connected,bujtas2019connected}, fractional \cite{bujtas2019fractional}, and disjoint \cite{bujtas2016disjoint}
domination games, the Z-, L-, LL-games \cite{brevsar2019variety}, the Maker-Breaker domination game \cite{duchene2020maker}, and the Maker-Breaker total domination game \cite{gledel2020maker}. For games played on hypergraphs, there is the aforementioned transversal game \cite{bujtas2016transversal,bujtas2017bounds}, but also the domination game itself can also be played on hypergraphs \cite{bujtas2019domination}. All these variants have been inspired by research on the domination game, but before the latter was first studied in 2010, there had already been research on further optimization games that we want to reference due to their close similarity to the domination game. These are the oriented domination game \cite{alon2002game} the competition-independence game, and the enclaveless game \cite{phillips2001introduction,phillips2002graph}.
We refer the reader to \cite{brevsar2021domination} for a general survey about domination games. 

\section*{Acknowledgement}

The authors would like to thank Vojtěch Dvořák for providing us with interesting insights into the topic and the suggestion of using linear programming. We would also like to thank Jan Petr for pointing out \Cref{jans-lemma} early on, which turned out to be vital for the proof, Béla Bollobás for his valuable comments, and Douglas Rall for pointing out an inaccuracy in \Cref{SectionPhase1} in an earlier version of this paper.

\bibliographystyle{abbrvnat}  
\renewcommand{\bibname}{Bibliography}
\bibliography{bibliography}

\end{document}